\documentclass[12pt]{article}
\usepackage{array,amsmath,amssymb}
\renewcommand{\baselinestretch}{1.13}           

\begin{document}


\input epsf.tex
\def\sp{\bigskip}
\def\ti{\\ \hglue \the \parindent}
\def\ce#1{\LP\medskip\centerline{#1}\medskip}
\def\LP{\par\noindent}

\newtheorem{theorem}{Theorem}[section]
\newtheorem{lemma}[theorem]{Lemma}
\newtheorem{corollary}[theorem]{Corollary}
\newtheorem{proposition}[theorem]{Proposition}
\newtheorem{definition}[theorem]{Definition}
\newtheorem{conjecture}[theorem]{Conjecture}
\newtheorem{remark}[theorem]{Remark}
\newtheorem{example}[theorem]{Example}
\newtheorem{question}[theorem]{Question}
\def\TH#1{\begin{theorem}{\rm #1}\end{theorem}}
\def\LM#1{\begin{lemma}{\rm #1}\end{lemma}}
\def\CO#1{\begin{corollary}{\rm #1}\end{corollary}}
\def\PR#1{\begin{proposition}{\rm #1}\end{proposition}}
\def\DF#1{\begin{definition}{\rm #1}\end{definition}}
\def\CJ#1{\begin{conjecture}{\rm #1}\end{conjecture}}
\def\RK#1{\begin{remark}{\rm #1}\end{remark}}
\def\PF{\LP{\bf Proof.}}
\def\qed{\hfill\rule{2mm}{2mm}\par\bigskip}
\newenvironment{proof}{\PF}{\qed}

\def\al{\alpha} \def\be{\beta}  \def\ga{\gamma} \def\dlt{\delta}
\def\eps{\epsilon} \def\th{\theta}  \def\ka{\kappa} \def\lmb{\lambda}
\def\sg{\sigma} \def\om{\omega} \def\GA{\Gamma} \def\DLT{\Delta}
\def\LMB{\Lambda}  \def\SG{\Sigma}  \def\OM{\Omega}

\def\symd{{\scriptstyle\triangle}}    \def\bs{\backslash}
\def\bu{\bullet} \def\isom{\cong}  \def\nul{\varnothing}

\def\NN{{\Bbb N}} \def\RR{{\Bbb R}} \def\ZZ{{\Bbb Z}} \def\FF{{\Bbb F}}

\def\st{\colon\,}   

\def\esub{\subseteq}    \def\esup{\supseteq} \def\nosub{\not\esub}
\def\adj{\leftrightarrow} \def\nadj{\not\leftrightarrow}
\def\Eqv{\Leftrightarrow} \def\Imp{\Rightarrow} \def\implied{\Leftarrow}
\def\da{\downarrow}  \def\ua{\uparrow}

\def\({\left(}  \def\){\right)}
\def\[{\left[}  \def\]{\right]}
\def\la{\langle}        \def\ra{\rangle}
\def\lb{\left[}         \def\rb{\right]}
\def\lc{\lceil}         \def\rc{\rceil}
\def\lf{\lfloor}        \def\rf{\rfloor}

\def\BR#1#2{{{#1}\brack {#2}}}
\def\ANGLE#1#2{{#1 \atopwithdelims <> #2}}
\def\CH#1#2{\binom{#1}{#2}} \def\MULT#1#2#3{\binom{#1}{#2,\ldots,#3}}
\def\FR#1#2{\frac{#1}{#2}}
\def\pile#1#2{{#1 \atop #2}}
\def\FL#1{\left\lfloor{#1}\right\rfloor} \def\FFR#1#2{\FL{\frac{#1}{#2}}}
\def\CL#1{\left\lceil{#1}\right\rceil}   \def\CFR#1#2{\CL{\frac{#1}{#2}}}

\def\SE#1#2#3{\sum_{#1=#2}^{#3}}  \def\SGE#1#2{\sum_{#1\ge#2}}
\def\PE#1#2#3{\prod_{#1=#2}^{#3}} \def\PGE#1#2{\prod_{#1\ge#2}}
\def\UE#1#2#3{\bigcup_{#1=#2}^{#3}}

\def\SM#1#2{\sum_{#1\in #2}}    \def\PM#1#2{\prod_{#1\in #2}}
\def\UM#1#2{\bigcup_{#1\in #2}} \def\IM#1#2{\bigcap_{#1\in #2}}
\def\SC#1#2{\sum_{#1\esub #2}}  \def\PC#1#2{\prod_{#1\esub #2}}

\def\VEC#1#2#3{#1_{#2},\ldots,#1_{#3}}
\def\VECOP#1#2#3#4{#1_{#2}#4\cdots #4 #1_{#3}}
\def\LINF#1{\lim_{#1\to\infty}}
\def\LIM#1#2{\lim_{#1\to#2}}
\def\MAP#1#2#3{#1\colon\,#2\to#3}

\def\C#1{\left | #1 \right |}    


\def\ph{\widehat{p}}
\def\pecn{p}
\def\pec{p}
\def\specn{\ph}
\def\ov#1{{\overline{#1}}}
\def\diam{{\rm diam\,}}
\def\twolgn{2^{\CL{\lg n}}}

\title{Parity Edge-Coloring of Graphs}

\author{David P. Bunde\thanks{Department of Computer Science,
University of Illinois, Urbana, IL, bunde@cs.uiuc.edu.  Partially supported by
NSF grant CCR 0093348.}\,,
Kevin Milans\thanks{Department of Computer Science,
University of Illinois, Urbana IL, milans@cs.uiuc.edu}\,,
Douglas B. West\thanks{Department of Mathematics,
University of Illinois, Urbana, IL, west@math.uiuc.edu.  Work supported in part
by the NSA under Award No.~MDA904-03-1-0037.}\,,
Hehui Wu\thanks{Department of Mathematics,
University of Illinois, Urbana, IL.}
}
\date{}
\maketitle

\begin{abstract}
A {\it parity walk} in an edge-coloring of a graph is a walk along which each
color is used an even number of times.  We introduce two parameters.  Let
$\pecn(G)$ be the least number of colors in an edge-coloring of $G$ having no
parity path (a {\it parity edge-coloring}).  Let $\specn(G)$ be the least
number of colors in an edge-coloring of $G$ in which every parity walk is
closed (a {\it strong parity edge-coloring}).  Always
$\specn(G)\ge \pecn(G)\ge \chi'(G)$.

The main result is that $\specn(K_n)=2^{\CL{\lg n}}-1$ for all $n$.
Furthermore, the optimal coloring for $K_n$ is unique when $n$ is a power of 2
and completely described for all $n$.  Also $\pecn(K_n)=\specn(K_n)$ when
$n\le16$.  The main result strengthens a special case of a result of Daykin
and Lov\'asz on Boolean functions.

A connected graph $G$ lies in the hypercube $Q_k$ if and only if $G$ has a
parity $k$-edge-coloring in which every cycle is a parity walk.  Hence
$\pecn(G)\ge\CL{\lg n(G)}$, with equality for paths and even cycles.  When $n$
is odd, $\pecn(C_n)=\specn(C_n)=1+\CL{\lg n}$.  Also,
$\pecn(K_{2,n})=\specn(K_{2,n})$, with value $n$ when $n$ is even and $n+1$
when $n$ is odd.  In general, $\specn(K_{m,n})\le m'\CL{n/m'}$, where
$m'=2^{\CL{\lg m}}$.

Let $\pecn_r(G)$ be the least number of colors needed to assign $r$ colors to
each edge of $G$ so that every choice of a color from the list assigned to each
edge yields a parity edge-coloring.  Trivially, $\pecn_r(G)\le r\pecn(G)$;
we prove that equality holds for paths.
\end{abstract}

\section{Introduction}

Our work began by studying which graphs embed in the hypercube $Q_k$, the graph
with vertex set $\{0,1\}^k$ in which vertices are adjacent when they differ in
exactly one coordinate.  Mitas and Reuter~\cite{MR} motivated that question by
observing that the hypercube is a common architecture for parallel computing.
Coloring each edge with the position of the bit in which its endpoints differ
yields two necessary conditions for the coloring inherited by a subgraph $G$:
\looseness -1

\medskip
1) every cycle uses each color an even number of times,

\smallskip
2) every path uses some color an odd number of times.

\medskip

\noindent
The characterization that having a $k$-edge-coloring satisfying conditions (1)
and (2) is also sufficient for a connected graph $G$ to be a subgraph of $Q_k$
was proved as early as 1972, by Havel and Mov\'arek~\cite{HM}.  The problem was
studied as early as 1953 by Shapiro~\cite{Sha}.

Define the {\it usage} of a color on a walk to be the parity of the number of
times it appears along the walk.  A {\it parity walk} is a walk in which the
usage of every color is even.  Condition (1) for an edge-coloring
states that every cycle is a parity walk, and a stronger version of (2) is
the statement that every parity walk is closed.

In general, define a {\it parity edge-coloring} to be an edge-coloring having
no parity path.  Although some graphs do not embed in any hypercube, using
distinct colors on the edges produces a parity edge-coloring for any graph.
Hence we introduce the {\it parity edge-chromatic number} $\pecn(G)$, defined
to be the minimum number of colors in a parity edge-coloring of $G$.  Paths of
length 2 guarantee that every parity edge-coloring is a proper edge-coloring,
and hence $\pecn(G)\ge \chi'(G)$, where $\chi'(G)$ denotes the edge-chromatic
number.

A more restricted edge-coloring notion has more robust algebraic properties.
Define a {\it strong parity edge-coloring (spec)} to be an edge-coloring in
which every parity walk is closed.  Again using distinct colors produces such a
coloring, and we introduce the {\it strong parity edge-chromatic number}
$\specn(G)$, defined to be the minimum number of colors in a spec.  Since a
path is an open walk (that is, the endpoints are distinct), a spec has no
parity path.  Hence every spec is a parity edge-coloring, and
$\specn(G)\ge\pecn(G)$ for every graph $G$.

The characterization of subgraphs of $Q_k$ yields $\pecn(G)\ge \CL{\lg n(G)}$
when $G$ is connected, with equality for a path or even cycle (here $n(G)$
denotes $|V(G)|$).  When $n$ is odd, $\pecn(C_n)=\specn(C_n)=1+\CL{\lg n}$.
Also $\pecn(K_{2,n})=\specn(K_{2,n})$, with value $n$ when $n$ is even and
$n+1$ when $n$ is odd.  In these examples,
$\pecn(G)=\specn(G)$; we also give examples where equality does not hold.

Our main result concerns complete graphs: $\specn(K_n)=2^{\CL{\lg n}}-1$.
To motivate our focus on complete graphs, we note that this result strengthens
a special case of an old result in extremal set theory.  Daykin and Lov\'asz
\cite{DL} proved that if $S$ is a family of finite sets $B$ is a nontrivial
Boolean function, then $\{B(u,v)\st u,v\in S\}$ has size at least $\C S$.
Marica and Sch\"onheim \cite{MS} earlier proved the special case where $B$ is
set difference.  Our result strengthens the conclusion in the case where $B$ is
symmetric difference.

\begin{theorem}\label{setappl}
If $S$ is a family of $n$ finite sets, and $\oplus$ denotes symmetric
difference, then the size of $\{u\oplus v\st u,v\in S\}$ is at least
$2^{\CL{\lg n}}$.  Equality holds for every family of $n$ subsets of a set
of size $\CL{\lg n}$.
\end{theorem}
\begin{proof}
View each member of $S$ as a vertex of $K_n$.  Color $E(K_n)$ by assigning the
symmetric difference $u\oplus v$ (that is, the binary sum of the incidence
vectors) to the edge $uv$.  Consider a parity walk starting from vertex $u$.
As an edge with a particular color is traversed, each element in the name of
that color is added to or deleted from the name of the current vertex to reach
the name of the next vertex.  Since each color is used an even number of times,
the last vertex is the same as $u$.  That is, every parity walk is closed, and
the coloring is a spec.

By our main result, at least $2^{\CL{\lg n}}-1$ colors are used, and these
colors all denote nonempty sets.  Also $u\oplus u=\nul$ for all $u\in S$.
\end{proof}

Daykin and Lov\'asz noted that if $\C S=n$, then at least $n-1$ distinct
nonempty sets arise as symmetric differences of members of $S$, and if $n$ is
not a power of 2, then there are more.  Our result gives the optimal bound for
all $n$.  They also hinted that having only $n-1$ symmetric differences of
distinct members may require $S$ to have a special structure, presaging the
uniqueness that we prove for the optimal spec of $K_{2^k}$.

The ideas needed for the main result are algebraic.  Relative to a given
$k$-edge-coloring, the {\it parity vector} $\pi(W)$ of a walk $W$ is the binary
$k$-tuple whose $i$th bit is the usage of color $i$ along $W$ (expressed as 0
for even and 1 for odd).  We study the binary vector space $L_f$ consisting of
the parity vectors of closed walks relative to a given spec $f$.  When $f$ is a
spec of $K_n$, all nonzero vectors in $L_f$ have at least two 1s.

We use these properties to show that if some color in an optimal spec of $K_n$
is not used on a perfect matching, then $\specn(K_{n+1})=\specn(K_n)$.  On the
other hand, if every color class is a perfect matching, then $n$ is a power of
2 and the coloring is isomorphic to the {\it canonical edge-coloring}, where
the vertices are named by distinct binary $(\lg n)$-tuples and the color on the
edge $uv$ is $u+v$, using binary vector addition.  When $n$ is not a power of
2, every optimal spec of $K_n$ is obtained by deleting vertices from a
canonically colored $K_{2^{\CL{\lg n}}}$.

The complete bipartite graph $K_{n,n}$ behaves like $K_n$ in that
$\pecn(K_{n,n})=\specn(K_{n,n})=\chi'(K_{n,n})= n$ when $n=2^k$.
Also, $\specn(K_{n,n})\le \specn(K_n)+1$ for all $n$; we conjecture that
equality holds.  We show that $\specn(K_{m,n})\le m'\CL{n/m'}$, where $m\le n$
and $m'=2^{\CL{\lg m}}$.

We have computed $\specn(K_n)$, but we do not know whether
$\pecn(K_n)=\specn(K_n)$ for all $n$ (and similarly for $K_{n,n}$).  As a
possible tool, we generalize the notion of parity edge-coloring.
A {\it parity $r$-set edge-coloring} assigns $r$ colors to each edge so that
every selection of one color from the set at each edge yields a parity
edge-coloring.  Let $p_r(G)$ be the minimum total number of colors used.
Always $\pecn_r(G)\le r\pecn(G)$, and we prove equality for paths.  Proving
$\pecn_2(K_n)=2\pecn(K_n)$ could be a step toward proving
$\pecn(K_n)=2^{\CL{\lg n}}-1$.

In the penultimate section of the paper, we describe related edge-coloring
problems with other constraints for the usage of colors on paths, and we
distinguish $\pecn(G)$ from those.  The final section poses many open questions.

\section{Elementary Properties and Examples}

First we state formally some elementary observations from the Introduction.

\begin{remark}
For every graph $G$, $\specn(G)\ge \pecn(G)\ge \chi'(G)$, and the parameters
$\specn$ and $\pecn$ are monotone under the subgraph relation.
\end{remark}
\begin{proof}
As noted earlier, $\pecn(G)\ge \chi'(G)$ by considering paths of length 2, and
$\specn(G)\ge \pecn(G)$ since closed walks are not paths.  When $H\esub G$, a
parity edge-coloring or spec of $G$ restricts to an edge-coloring of that type
on $H$, since every parity walk in the restriction to $H$ is a parity walk in
the original coloring of $G$.
\end{proof}

When $G$ is a forest, every parity edge-coloring is also a spec, so
$\pecn(G)=\specn(G)$.  We have observed that the edge-coloring of the hypercube
by coordinates shows that $\specn(Q_k)=\pecn(Q_k)=k$.  Hence $\specn(G)\le k$
for every subgraph $G$ of $Q_k$.  For trees, this will also be sufficient.

Recall that given a $k$-edge-coloring $f$ and a walk $W$, $\pi(W)$ denotes the
parity vector of $W$, recording the usage of each color as 0 or 1.  When walks
$W$ and $W'$ are concatenated, the parity vector of the concatenation is the
vector binary sum $\pi(W)+\pi(W')$.  The {\it weight} of a vector is the
number of nonzero positions.

\begin{theorem}\label{tree}
A tree $T$ embeds in the $k$-dimensional hypercube $Q_k$ if and only if
$\pecn(T)\le k$.
\end{theorem}
\begin{proof}
We have observed necessity.  Conversely, let $f$ be a parity $k$-edge-coloring
of $T$ (there may be unused colors if $\pecn(T)<k$).  Fix a root vertex $r$ in
$T$.  Define $\MAP\phi{V(T)}{V(Q_k)}$ by setting $\phi(v)=\pi(W)$, where $W$ is
the $r,v$-path in $T$.

When $uv\in E(T)$, the $r,u$-path and $r,v$-path in $T$ differ in one edge,
so $\phi(u)$ and $\phi(v)$ are adjacent in $Q_k$.  It remains only to check
that $\phi$ is injective.  The parity vector for the $u,v$-path $P$ in $T$ is
$\phi(u)+\phi(v)$, since summing the $r,u$-path and $r,v$-path cancels the
portion from $r$ to $P$.  Since $f$ is a parity edge-coloring, $\phi(P)$ is
nonzero, and hence $\phi(u)\ne\phi(v)$.
\end{proof}

When $k$ is part of the input, recognizing subgraphs of $Q_k$ is NP-complete
\cite{KVC}, and this remains true when the input is restricted to trees
\cite{embedding-hardness}.  Therefore, computing $\pecn(G)$ or $\specn(G)$ is
NP-hard even when $G$ is a tree.  Perhaps there is a polynomial-time algorithm
for trees with bounded degree or bounded diameter.

Havel~\cite{Havel} proposed studying the spanning trees of $Q_k$, and many
papers quickly followed; Havel~\cite{Havel2} presents a survey.  It is
necessary that the tree be {\it equitable} (viewed as a bipartite graph, its
partite sets have the same size), but generally this is not sufficient.
Kobeissi and Mollard \cite{KM} proved that it is sufficient for a
{\it double-starlike tree}, which is a subdivision of a double star in which
the central edge is not subdivided.  This result strengthened a string of
earlier results (such as \cite{Neb2}) on the special case of subdivisions of
stars.  Equitability and order $2^k$ are also sufficient for various classes
of caterpillars (see \cite{DHLM,HavLie}).

Next we consider arbitrary subgraphs of $Q_k$.  As noted earlier, Havel and
Mov\'arek\cite{HM} proved that having a $k$-edge-coloring with properties (1)
and (2) of the Introduction is necessary and sufficient for a graph to be a
subgraph of $Q_k$.  (They also proved statements equivalent to
Theorem~\ref{tree} and Corollary~\ref{pathcyc}.)  Their proof is essentially
the same as ours, though our organization is different in order to motivate the
parameters we have defined.  We phrase the condition using parity edge-coloring
and parity walks and express the result as a corollary of Theorem~\ref{tree} to
motivate our later arguments.

\begin{corollary}
A graph $G$ is a subgraph of $Q_k$ if and only if $G$ has a parity
$k$-edge-coloring in which every cycle is a parity walk.
\end{corollary}
\begin{proof}
We have observed necessity.  For sufficiency, choose a spanning tree $T$.
Since $\pecn(T)\le \pecn(G)\le k$, Theorem~\ref{tree} implies that $T\esub Q_k$.
Map $T$ into $Q_k$ using $\phi$ as defined in the proof of Theorem~\ref{tree}.
For each $xy\in E(G)-E(T)$, the cycle formed by adding $xy$ to $T$ is given to
be a parity walk.  Hence the $x,y$-path in $T$ has parity vector with weight 1.
This makes $\phi(x)$ and $\phi(y)$ adjacent in $Q_k$, as desired.
\end{proof}

Mitas and Reuter~\cite{MR} later gave a much lengthier proof motivated by
applying analogous methods to study subdiagrams of the subset lattice.
They also characterized the graphs occurring as induced subgraphs of $Q_k$
as those having a $k$-edge-coloring satisfying properties (1) and (2) and (3),
where property (3) essentially states that that if the parity vector of a walk
$W$ has weight 1, then the endpoints of $W$ are adjacent.

Spanning trees yield a general lower bound on $\pecn(G)$, which holds with
equality for paths, even cycles, and connected spanning subgraphs of $Q_k$.

\begin{corollary}\label{connected}
If $G$ is connected, then $\pecn(G) \geq \CL{\lg n(G)}$.
\end{corollary}
\begin{proof}
If $T$ is a spanning tree of $G$, then $\pecn(G) \geq \pecn(T)$.  Since
$T$ embeds in the hypercube of dimension $\pecn(T)$,
we have $n(G)=n(T) \le 2^{\pecn(T)}\le 2^{\pecn(G)}$.
\end{proof}

\begin{corollary}\label{pathcyc}
For all $n$, $\pecn(P_n) = \specn(P_n) = \CL{\lg n}$.
For even $n$, $\pecn(C_n) = \specn(C_n) = \CL{\lg n}$.
\end{corollary}
\begin{proof}
The lower bounds follow from Corollary~\ref{connected}.  The upper bounds
hold because $Q_k$ contains cycles of all even lengths up to $2^k$.
\end{proof}

A result equivalent to $\pecn(P_n)=\specn(P_n)=\CL{\lg n}$ appears in
\cite{HM} (without defining either parameter).  When $n$ is odd, $C_n$ needs an
extra color beyond $\CL{\lg n}$.  To prove this, we begin with simple
observations about adding an edge.

\begin{lemma}\label{addedge}
(a) If $e$ is an edge in a graph $G$, then $\pecn(G)\le \pecn(G-e)+1$.\\
(b) If also $G-e$ is connected, then $\specn(G)\le\specn(G-e)+1$.
\end{lemma}
\begin{proof}
(a) Put an optimal parity edge-coloring on $G-e$ and add a new color on $e$.
There is no parity path avoiding $e$, and any path through $e$ uses the
new color exactly once.

(b) Put an optimal spec on $G-e$ and add a new color on $e$.
Let $P$ be a $u,v$-path in $G-e$, where $u$ and $v$ are the endpoints of $e$.
Suppose that there is an open parity walk $W$.  Note that $W$ traverses $e$
an even number of times, since no other edge has the same color as $e$.
Form $W'$ by replacing each traversal of $e$ by $P$ or its reverse, depending
on the direction of traversal of $e$.  Every edge is used with the same parity
in $W'$ and $W$, and the endpoints are unchanged, so $W'$ is an open parity
walk in $G-e$.  This is a contradiction.
\end{proof}

Lemma~\ref{addedge}(b) does not hold when $G-e$ is disconnected (see
Example~\ref{phnep}).

\begin{theorem}\label{oddcycle}
If $n$ is odd, then $\pecn(C_n) = \specn(C_n) = \CL{\lg n} + 1$.
\end{theorem}
\begin{proof}
Lemma~\ref{addedge}(b) yields the upper bound, since $\specn(P_n)=\CL{\lg n}$.

For the lower bound, we show first that $\specn(C_n) = \pecn(C_n)$ (this and
Lemma~\ref{addedge}(a) yield an alternative proof of the upper bound).  Let $W$
be an open walk, and let $W'$ be the subgraph formed by the edges with odd
usage in $W$.  The sum of the usage by $W$ of edges incident to a vertex $x$ is
odd if and only if $x$ is an endpoint of $W$.  Hence $W'$ has odd degree
precisely at the endpoints of $W$.  Within $C_n$, this requires $W'$ to be
a path $P$ joining the endpoints of $W$.  Under a parity edge-coloring $f$,
some color has odd usage along $P$, and this color has odd usage in $W$.  Hence
$f$ has no open parity walk, and every parity edge-coloring is a spec.

It now suffices to show that $\specn(C_n) \ge \pecn(P_{2n})$.  Given a spec $f$
of $C_n$, we form a parity edge-coloring $g$ of $P_{2n}$ with the same number
of colors.  Let $\VEC v1n$ be the vertices of $C_n$ in order, and let
$\VEC u1n,\VEC w1n$ be the vertices of $P_{2n}$ in order.  Define $g$ by
letting $g(u_iu_{i+1}) = g(w_iw_{i+1}) = f(v_iv_{i+1})$ for $1\le i\le n-1$ and
letting $g(u_nw_1) = f(v_nv_1)$.

Each path in $P_{2n}$ corresponds to an open walk in $C_n$ or to one trip
around the cycle.  There is no parity path of the first type, since $f$ is a
spec.  There is none of the second type, since $C_n$ has odd length.
\end{proof}

The ``unrolling'' technique of Theorem~\ref{oddcycle} leads to an example $G$
with $\specn(G) > \pecn(G)$, which easily extends to generate infinite families.

\begin{example}\label{phnep}
\rm
Form a graph $G$ by identifying a vertex of $K_3$ with an endpoint of $P_8$.
Since $\pecn(K_3)=\pecn(P_7)=3$, adding the connecting edge yields
$\pecn(G)\le 4$ (see Lemma~\ref{addedge}(a)).

We claim that $\specn(G)\ge \pecn(P_{18})=5$.  We copy a spec $f$ of $G$ onto
$P_{18}$ with the path edges doubled.  Beginning with the vertex of degree 1
in $G$, walk down the path, once around the triangle, and back up the path.
This walk has length 17; copy the colors of its edges in order to the edges
of $P_{18}$ in order to form an edge-coloring $g$ of $P_{18}$.

Each path in $P_{18}$ corresponds to an open walk in $G$ or a closed walk
that traverses the triangle once.  There is no parity path of the first type,
since $f$ is a spec.  There is none of the second type, since such a closed
walk has odd length.  This proves the claim.

Since $\specn(K_3)=\specn(P_7)=3$, this graph $G$ also shows that adding
an edge can change $\specn$ by more than 1 when $G$ is disconnected.
\end{example}

We know of no bipartite graph $G$ with $\specn(G)>\pecn(G)$.
Nevertheless, it is not true that every optimal parity edge-coloring of
a bipartite graph is a spec.

\begin{example}
\rm
Let $G$ be the graph obtained from $C_6$ by adding two pendant edges at
one vertex.  Let $W$ be the spanning walk that starts at one pendant vertex,
traverses the cycle, and ends at the other pendant vertex.  Let $f$ be the
4-edge-coloring that colors the edges of $W$ in order as $a,b,a,c,b,d,c,d$.
Although $f$ is an optimal parity edge-coloring ($\Delta(G)=4$), it uses
each color twice on the open walk $W$, so it is not a spec.  Changing the
edge of color $d$ on the cycle to color $a$ yields a strong parity
$4$-edge-coloring.
\end{example}

\section{Complete Graphs and Linear Algebra}\label{linalg}

In this section we use linear algebra to prove our main result, determining
$\specn(K_n)$ for all $n$.  We begin with a construction when $n$ is a power of
$2$, for which we recall a definition from the Introduction.

\begin{definition}
\rm
When $n=2^k$, the {\em canonical coloring} of $K_n$ is the edge-coloring $f$
defined by $f(uv)=u+v$, where $V(K_n)=\FF_2^k$ and addition is binary vector
addition.
\end{definition}

\begin{lemma}\label{canon}
If $n=2^k$, then $\specn(K_n)=\pecn(K_n)=\chi'(K_n)=n-1$.
\end{lemma}
\begin{proof}
The canonical coloring $f$ uses $n-1$ colors (the color $0^k$ is not used).
We show that $f$ is a spec.  When $W$ is an open walk, its ends differ in some
bit $i$.  The total usage of colors flipping bit $i$ along $W$ is odd, and
hence some color has odd usage on $W$.
\end{proof}

Since every complete graph is a subgraph of the next larger complete graph,
we obtain $\specn(K_n) \leq 2^{\CL{\lg n}} - 1$.  We will show that this
upper bound is exact.  The main idea is that we will be able to introduce
an additional vertex without needing additional colors until a power of 2 is
reached.  At that point, Theorem~\ref{canonical} will apply.

\begin{definition}\label{4condef}
\rm
An edge-coloring $f$ of $G$ satisfies the {\em 4-constraint} if whenever
$f(uv)=f(xy)$ and $vx\in E(G)$, also $uy\in E(G)$ and $f(uy)=f(vx)$.
\end{definition}

\begin{lemma}\label{4constraint}
If $f$ is a parity edge-coloring in which every color class is a perfect
matching, then $f$ satisfies the 4-constraint.
\end{lemma}
\begin{proof}
Otherwise, given $f(uv)=f(xy)$, the edge of color $f(vx)$ incident to $u$ forms
a parity path of length 4 with $uv$, $vx$, and $xy$.
\end{proof}

\begin{theorem}\label{canonical}
If $f$ is a parity edge-coloring of $K_n$ in which every color class is a
perfect matching, then $f$ is a canonical coloring and $n$ is a power of $2$.
\end{theorem}
\begin{proof}
Every edge is a canonically colored copy of $K_2$.  Let $R$ be a largest
vertex set on which $f$ restricts to a canonical coloring, so $|R|=2^{j-1}$ for
some $j$.  We are given a bijection $\phi$ from $R$ to $\FF_2^{j-1}$ under
which $f$ is the canonical coloring.

Since $f$ is canonical, every color used within $R$ by $f$ pairs the vertices
of $R$.  Let $c$ be a color not used within $R$; since $c$ is used on a
perfect matching, $c$ matches $R$ to some set $U$.  Let $R'=R\cup U$.
Define $\MAP{\phi'}{R'}{\FF_2^j}$ as follows: for $x\in R$, obtain $\phi'(x)$
by appending $0$ to $\phi(x)$; for $x\in U$ obtain $\phi'(x)$ by appending $1$
to $\phi(x')$, where $x'$ is the neighbor of $x$ in color $c$.
Within $R'$, we henceforth refer to the vertices by their names under $\phi'$.

By Lemma~\ref{4constraint}, the 4-constraint holds for $f$.  The 4-constraint
copies the coloring from the edges within $R$ to the edges within $U$.  To see
this, consider $x',y'\in U$ arising from $x,y\in R$, with $f(xx')=f(yy')=c$.
Now $f(x'y')=f(xy)=x+y=x'+y'$, using the 4-constraint, the fact that $f$ is
canonical on $R$, and the definition of $\phi'$.  Hence $f$ is canonical within
$U$.

Finally, let $u$ be the name of the color on the edge $0^j u$, for $u\in U$.
For any $v\in R$, let $w=u+v$; note that $w\in U$.  Both $0^jv$ and $uw$
have color $v$, since $f$ is canonical within $R$ and within $U$.  By applying
the 4-constraint to $\{v0^j,0^jw,wu\}$, we conclude that $f(uv)=f(0^jw)=w$.
Since $w=u+v$, this completes the proof that $f$ is canonical on $R'$.
\end{proof}

Now we begin the algebraic observations needed to prove the main result.

\begin{lemma}\label{vspace}
For an edge-coloring $f$ of a connected graph $G$, the set $L$ of parity
vectors of closed walks is a vector space under binary vector addition.
\end{lemma}
\begin{proof}
Traversing an edge twice yields a closed walk with parity vector zero.
It thus suffices to show that $L$ is closed under addition over $\FF_2$.
Given a $u,u$-walk $W$ and a $v,v$-walk $W'$, let $P$ be a $u,v$-path in $G$,
and let $\ov{P}$ be its reverse.  Now following $W,P,W',\ov{P}$ in succession
yields a $u,u$-walk with parity vector $\pi(W)+\pi(W')$.
\end{proof}

\begin{definition}\label{parityspace}
\rm
For an edge-coloring $f$, the {\em parity space} $L_f$ is the vector space
of parity vectors of closed walks under $f$.  Let $w(L)$ denote the minimum
weight of a nonzero vector in a binary vector space $L$.
\end{definition}

\begin{lemma}\label{weight2}
If an edge-coloring $f$ of a graph $G$ is a spec, then $w(L_f)\ge2$.
The converse holds when $G=K_n$.
\end{lemma}
\begin{proof}
If the parity vector of a closed walk $W$ has weight 1, then one color has odd
usage in $W$ (say on edge $e$).  Now $W-e$ is an open parity walk, and $f$
is not a spec.

If $f$ is not a spec, then there is an open parity walk $W'$.  In $K_n$, the
ends of $W'$ are adjacent, and adding that edge yields a closed walk whose
parity vector has weight 1.
\end{proof}

\begin{lemma}\label{3colors}
For colors $a$ and $b$ in an optimal spec $f$ of $K_n$, there is some closed
walk $W$ on which the colors having odd usage are $a$, $b$, and one other.
\end{lemma}
\begin{proof}
We use Lemma~\ref{weight2} repeatedly.  Since $f$ is optimal, merging the
colors $a$ and $b$ into a single color $a'$ yields an edge-coloring $f'$ that
is not a spec.  Hence under $f'$ there is a closed walk $W$ on which $f'$ has
odd usage for only one color $c$.  Also $c\ne a'$, since otherwise $f$ has odd
usage on $W$ for only $a$ or $b$.  With $c\ne a'$ and the fact that $f$ has odd
usage for at least two colors on $W$, both $a$ and $b$ also have odd usage on
$W$, and $W$ is the desired walk.
\end{proof}

The same idea as in Lemma~\ref{3colors} shows that $w(L_f)\ge3$ when $f$ is an
optimal spec of $K_n$, but we do not need this observation.  We note, however,
that the condition $w(L_f)\ge3$ is the condition for $L_f$ to be the set of
codewords for a 1-error-correcting code.  Indeed, when $n=2^k$ and $f$ is the
canonical coloring, $L_f$ is a perfect 1-error-correcting code of length $n-1$.

A {\em dominating vertex} in a graph is a vertex adjacent to all others.

\begin{lemma}\label{domvert}
If $f$ is an edge-coloring of a graph $G$ with a dominating vertex $v$,
then $L_f$ is the span of the parity vectors of triangles containing $v$.
\end{lemma}
\begin{proof}
By definition, the span is contained in $L_f$.  Conversely, consider any
$\pi(W)\in L_f$.  Let $S$ be the set of edges with odd usage in $W$, and let
$H$ be the spanning subgraph of $G$ with edge set $S$.  Since the total
usage at each vertex of $W$ is even, $H$ is an even subgraph of $G$.
Hence $H$ decomposes into cycles, which are closed walks, and $\pi(W)$ is
the sum of the parity vectors of these cycles.

It therefore suffices to show that $S$ is the set of edges that appear in an
odd number of the triangles formed by $v$ with edges of $H-v$.  Each edge of
$H-v$ is in one such triangle, so we need consider only edges involving $v$.
An edge $vw$ lies in an odd number of these triangles if and only if
$d_{H-v}(w)$ is odd, which occurs if and only if $w\in N_H(v)$, since
$d_H(w)$ is even.  By definition, $vw\in E(H)$ if and only if $vw$ has
odd usage in $W$ and hence lies in $S$.
\end{proof}

\begin{lemma}\label{nonmatch}
If an optimal spec $f$ of $K_n$ uses some color $a$ on less than a perfect
matching, then $\specn(K_{n+1})=\specn(K_n)$.
\end{lemma}
\begin{proof}
We view $K_{n+1}$ as arising from $K_n$ by adding a vertex $u$.  Let $v$ be a
vertex of $K_n$ at which $a$ does not appear.

We use $f$ to define $f'$ on $E(K_{n+1})$.  Let $f'$ agree with $f$ on
$E(K_n)$, and let $f'(uv)=a$.  To define $f'$ on each remaining edge $uw$,
first let $b=f(vw)$.  By Lemma~\ref{3colors}, there is a closed walk $W$ with
odd usage precisely for $a$ and $b$ and some third color $c$ under $f$.  Let
$f'(uw)=c$.

Note that $f'$ uses the same colors as $f$.  It remains only to show that $f'$
is a spec.  To do this we prove that $w(L_{f'})\ge2$, by showing that
$L_{f'}\esub L_f$.  By Lemma~\ref{domvert}, it suffices to show that
$\pi(T)\in L_f$ when $T$ is a triangle in $K_{n+1}$ containing $v$.

Triangles not containing $u$ lie in the original graph and have parity vectors
in $L_f$.  Hence we consider the triangle $T$ formed by $\{u,v,w\}$.
Now $\pi(T)=\pi(W)\in L_f$, where $W$ is the walk used to specify $f'(uw)$.
\end{proof}

\begin{theorem}\label{specKn}
$\specn(K_n)=2^{\CL{\lg n}}-1$.
\end{theorem}
\begin{proof}
If some color class in an optimal spec is not a perfect matching, then
$\specn(K_n)=\specn(K_{n+1})$, by Lemma~\ref{nonmatch}.  This vertex absorption
cannot stop before the number of vertices reaches a power of 2, because when
every color class is a perfect matching the coloring is canonical, by
Theorem~\ref{canonical}.  It cannot continue past $2^{\CL{\lg n}}$ vertices,
since then the maximum degree equals the number of colors.  Hence
$\specn(K_n)=\specn(K_{2^{\CL{\lg n}}})=2^{\CL{\lg n}}-1$.
\end{proof}

Although we do not know the complexity of recognizing parity edge-colorings,
our algebraic results settle that question for specs.

\begin{theorem}\label{checspec}
The problem of recognizing strong parity edge-colorings of graphs is solvable
in polynomial time.
\end{theorem}
\begin{proof}
Let $f$ be an edge-coloring of a graph $G$.  By treating each component of $G$
separately, we may assume that $G$ is connected.  Now we may also assume that
$G$ is a complete graph, since adding a missing edge and giving it a new color
does not change whether the coloring is a spec, by Lemma~\ref{addedge}(b).
Let $f'$ be the resulting coloring.

Next we obtain a spanning set for $L_{f'}$.  Every vertex of $K_n$ is a
dominating vertex, so by Lemma~\ref{domvert}, ${f'}_f$ is the span of the parity
vectors of triangles containing $v$.  There are $\CH{n-1}2$ such triangles, and
we obtain each parity vector in constant time.

We now form a matrix with these parity vectors as the columns.  By
Lemma~\ref{weight2}, it suffices to check whether any vector of weight 1
is in their span.  With Gaussian elimination, we can check all such vectors
in time polynomial in $n$.
\end{proof}

It is natural to wonder whether every edge-coloring of $K_n$ that satisfies
the 4-constraint is a spec or a parity edge-coloring.  The next example
shows that the answer is no.
Similarly, not every parity edge-coloring of $K_n$ is a spec.  Nevertheless,
proving that every {\em optimal} parity edge-coloring is a spec would prove the
conjecture that $\pecn(K_n)=\twolgn-1$.

\begin{example}
\rm
Color $E(K_n)$ as follows.  Let the colors on the edges of some spanning cycle
be $0,1,2,3,4,5,1,3,5,2,4$ in order.  No matter how the coloring is completed,
there is a parity path, so the coloring cannot be completed to a spec or even
a parity edge-coloring.  Nevertheless, we can complete it to satisfy the
4-constraint.

For $i\in\{1,2,3,4,5\}$, the vertices of the two edges with color $i$ induce
a copy of $K_4$.  The 4-constraint requires this $K_4$ to be colored with
three matchings; let the other two have colors $i+5$ and $i+10$.
Because no two colors are incident twice on the original cycle, these copies
of $K_4$ are pairwise edge-disjoint, so we can color them independently.
We have now colored $31$ edges.  We give the remaining 24 edges distinct new
colors, so there are no further requirements from the 4-constraint.

We have used 40 colors, although $\specn(K_{11})=15$.  It is possible that
the 4-constraint is sufficient for a spec of $K_n$ when the number of colors
is restricted near $\specn(K_n)$.
\end{example}

\section{Complete Bipartite Graphs}

Parity edge-colorings of complete bipartite graphs are related to those of
complete graphs.

\begin{proposition}\label{bicanon}
If $n=2^k$, then $\specn(K_{n,n}) = \pecn(K_{n,n}) = n$.
\end{proposition}
\begin{proof}
The lower bound is $\chi'(K_{n,n})$.  For the upper bound, we have the
``bicanonical coloring'' analogous to the canonical coloring of $K_n$.
Name the vertices of each partite set using the vectors from $\FF_2^n$,
and give $uv$ the color $u+v$.  As in Lemma~\ref{canon}, a parity walk
must start and end at the same label.  Since its length is even, it must
also start and end in the same partite set and hence at the same vertex.
\end{proof}

By the subgraph relation, $\specn(K_{n,n})\le 2^{\CL{\lg n}}$ for all $n$.
Toward the conjecture that equality holds, we offer the following.

\begin{proposition}\label{cliqbicliq}
If some optimal spec of $K_{n,n}$ uses a color on at least $n-1$ edges,
then $\specn(K_{n,n})=\specn(K_n)+1=2^{\CL{\lg n}}$.
If a color is used $n-r$ times, then $\specn(K_{n,n})\ge\twolgn-\CH r2$.
\end{proposition}
\begin{proof}
We prove the general statement.
Let $f$ be such a spec, and let $c$ be such a color.  Let $U$ be one partite
set, with vertices $\VEC u1n$.  Whenever color class $c$ is incident to at
least one of distinct vertices $u_i,u_j\in U$, let $P_{i,j}$ be a
$u_i,u_j$-path of length 2 in which one edge has color $c$ under $f$.  Choose
these so that $P_{j,i}$ is the reverse of $P_{i,j}$.  When $c$ appears at
neither $u_i$ nor $u_j$, leave $P_{i,j}$ undefined.

Let $G$ be the graph obtained from $K_n$ with vertex set $\VEC v1n$ by
deleting the edges $v_iv_j$ such that $P_{i,j}$ is undefined; there are
$\CH r2$ such edges.  Define a coloring $f'$ on $G$ by letting $f(v_iv_j)$ be
the color other than $c$ on $P_{i,j}$.

We claim that $f'$ is a spec.  Given a parity walk $W'$ under $f'$, define a
walk $W$ in $K_{n,n}$ as follows.  For each edge $v_iv_j$ in $W'$, follow
$P_{i,j}$.  By construction, the usage in $W$ of each color other than $c$ is
even.  Hence also the usage of $c$ is even.  Hence $W$ is a parity walk under
$f$ and therefore is closed.  Since $W$ starts and ends at the same vertex
$u_i\in U$, also $W'$ starts and ends at the same vertex $v_i$.

We have proved that every parity walk under $f'$ is closed, so $f'$ is a spec.
Hence $f'$ has at least $\specn(G)$ colors, and $f$ has at least one more.
By Lemma~\ref{addedge}(b) and Theorem~\ref{specKn},
$\specn(G)\ge\twolgn-1-\CH r2$, which completes the lower bound.

For the upper bound, Proposition~\ref{bicanon} shows that $\twolgn$ colors
suffice.
\end{proof}

\begin{corollary}\label{rbound}
$\specn(K_{n,n})\ge \max_r\min\{\twolgn-\CH r2,\FR{n^2}{n-r-1}\}$.
\end{corollary}
\begin{proof}
If $E(K_{n,n})$ has a spec with $s$ colors, where $s<\twolgn-\CH r2$, then
by Proposition~\ref{cliqbicliq} no color can be used at least $n-r$ times,
and hence ${n^2/s}\le n-r-1$.  Thus
$\specn(K_{n,n})\ge \min\{\twolgn-\CH r2,n^2/(n-r-1)\}$.
\end{proof}

With $r=1$, we conclude that $\specn(K_{n,n})\ge 2^k$ when $n>2^k-3-4/(n-2)$,
since then $n^2/(n-2)>2^k-1$.  Thus $\specn(K_{5,5})=8$, and
$\specn(K_{n,n})=16$ for $13\le n\le 16$.  Using $r=2$, we obtain
$14\le \specn(K_{9,9})\le 16$.

Lacking a direct proof that always $\specn(K_{n,n})=\specn(K_n)+1$, we note
that it may be possible to develop an algebraic proof of
$\specn(K_{n,n})=2^{\CL{\lg n}}$ analogous to that of Theorem~\ref{specKn}.
Doing so would {\em replace} much of Section~\ref{linalg} using the following
result.

\begin{proposition}
$\specn(K_n)\ge \specn(K_{n,n})-1$.
\end{proposition}
\begin{proof}
Let $f$ be a spec of $K_n$ with vertex set $\VEC u1n$.  Given $K_{n,n}$ with
partite sets $\VEC v1n$ and $\VEC w1n$, let $f'(v_iw_j)=f(u_iu_j)$ when
$i\ne j$, and give a single new color to all $v_iw_i$ with $1\le i\le n$.  A
parity walk $W'$ under $f'$ starts and ends in the same partite set.  Mapping
it back to $K_n$ (dropping the edges with the new color) yields an even parity
walk $W$ under $f$.  Hence $W$ starts and ends at vertices with the same index.
Since $K_{n,n}$ has only one vertex with each index in each partite set, $W'$
is closed.  Hence $f'$ is a spec of $K_{n,n}$.
\end{proof}

The conjectured value of $\specn(K_{n,n})$ would generalize
Theorem~\ref{setappl} as follows: If $S_1$ and $S_2$ both are families
consisting of $n$ finite sets, then
$\C{\{u\oplus v\st u\in S_1,v\in S_2\}}\ge \twolgn$.

Proposition~\ref{bicanon} shows that when $n$ is a power of $2$, the lower
bound of $\DLT(G)$ is optimal for strong parity edge-coloring of $K_{n,n}$.
We next enlarge the class of complete bipartite graphs where this bound is
optimal.

\begin{theorem}\label{Kmn}
If $m=2^k$ and $m$ divides $n$, then
$\pecn(K_{m,n})=\specn(K_{m,n})=\Delta(K_{m,n})=n$.
\end{theorem}
\begin{proof}
Let $r=n/m$ and $[r]=\{1,\ldots,r\}$.  Label the vertices in the small part
with $\FF_2^k$.  Label those in the large part with $\FF_2^k\times [r]$.
Color the edges with color set $\FF_2^k\times [r]$ by setting
$f(uv)=(u+v',j)$, where $v=(v',j)$.  In other words, we use $r$ edge-disjoint
copies of the bicanonical coloring on $r$ edge-disjoint copies of $K_{m,m}$.

We have used $n$ colors, so it suffices to show that $f$ is a spec.
Let $W$ be a parity walk under $f$.  Erasing the second coordinate maps
$W$ onto a walk $W'$ in $K_{m,m}$.  Furthermore, $W'$ is a parity walk, because
all edges in $W$ whose color has the form $(z,j)$ for any $j$ are mapped onto
edges with color $z$ under the bicanonical coloring of $K_{m,m}$, and there
are an even number of these for each $j$.  Hence $W'$ is closed.

Hence $W$ starts and ends at vertices labeled with the same element $u$ of
$\FF_2^k$, and they are in the same part since $W$ has even length.  If these
vertices are different copies of $u$ in the large partite set, then those
copies of $K_{m,m}$ have contributed an odd number of edges to $W$, so
for each of them some color confined to it has odd usage in $W$.  This
contradicts that $W$ is a parity walk.  Hence $W$ is closed, and $f$ is a spec.
\end{proof}

\begin{corollary}\label{morebicliq}
If $m\le n$ and $m' = 2^{\CL{\lg m}}$, then $\specn(K_{m,n})\le m'\CL{n/m'}$.
\end{corollary}
\begin{proof}
$K_{m,n}\esub K_{m',m'\CL{n/m'}}$.
\end{proof}

Corollary~\ref{morebicliq} provides examples of complete bipartite graphs where
the maximum degree bound is optimal even though the size of neither partite set
is a power of 2.  For example, $\specn(K_{3,12})=12$.  We use the corollary
next to compute the exact values when $m=2$.  We will apply the result for
$K_{2,3}$ in Theorem~\ref{K9}.

\begin{corollary}\label{K2n}
$\specn(K_{2,n})=\pecn(K_{2,n})$, with value $n$ for even $n$
and $n+1$ for odd $n$.
\end{corollary}
\begin{proof}
The upper bounds are given by Corollary~\ref{morebicliq}, where $m'=2$.

For the lower bound, since $\Delta(K_{2,n})=n$ for $n\ge2$, it suffices to show
that $n$ must be even when $f$ is a parity edge-coloring of $K_{2,n}$ with $n$
colors.  Let $\{x,x'\}$ be the partite set of size 2.  Each color appears at
both $x$ and $x'$.  If color $a$ appears on $xy$ and $x'y'$, then
$f(xy')=f(x'y)$, since otherwise the colors $a$ and $f(xy')$ form a parity path
of length 4.

Hence $y$ and $y'$ have the same pair of incident colors.  Making this
argument for each color partitions the vertices in the partite set of size $n$
into pairs.  Hence $n$ is even.
\end{proof}

The upper bound in Corollary~\ref{K2n} can also be proved using an augmentation
lemma.  If $f$ is a spec of a connected graph $G$, and $G'$ is formed from $G$
by adding new vertices $x$ and $y$ with common neighbors $u$ and $v$ in $G$
(and no other new edges), then the coloring $f'$ obtained from $f$ by adding
two new colors $a$ and $b$ alternating on the new 4-cycle is a spec of $G'$.
This yields $\specn(G')\le\specn(G)+2$.  Like Lemma~\ref{addedge}(b), this
statement fails for disconnected graphs.  Since we presently have no further
applications for this lemma, we omit the proof.

Before leaving the subject of strong parity edge-coloring, we observe that
every graph has an optimal spec satisfying a weaker form of the 4-constraint.

\begin{definition}
\rm
Given a spec $f$ of a graph $G$, let $m_i$ be the size of the $i$th color
class.  An optimal spec $f$ is {\em lex-optimal} if it lexicographically
maximizes the vector $(\VEC m1{\specn(G)})$ among all optimal specs.  An
edge-coloring $f$ satisfies the {\em weak 4-constraint} if $f(vx)=f(yu)$
whenever $f(uv)=f(xy)$ and $vx,yu\in E(G)$.
\end{definition}

The weak 4-constraint is weaker than the 4-constraint (Definition~\ref{4condef})
by not requiring 4-vertex paths with repeated colors to lie in 4-cycles.
For edge-colorings of complete graphs, they are equivalent.

\begin{proposition}\label{weak4}
For any graph $G$, a lex-optimal spec $f$ satisfies the weak 4-constraint.
\end{proposition}
\begin{proof}
If $f$ fails the weak 4-constraint, then there is a 4-cycle $C$ with vertices
$u,v,x,y$ such that $f(uv)=f(xy)=a$ but $b=f(vx)\ne f(yu)=c$.  By symmetry, we
may assume that $b$ is earlier than $c$ in the list of colors.  Obtain $f'$
from $f$ by changing the color of $yu$ from $c$ to $b$.

The vector of multiplicities for $f'$ lexicographically exceeds the vector
for $f$, so it suffices to show that $f'$ is a spec.  If not, then $f'$ yields
an open parity walk $W'$.  Since $f$ permits no such walk, $W'$ uses $yu$.
Replacing each traversal of $yu$ by a traversal of the rest of $C$ yields a
parity walk $W$ under $f$, since the usage of each color is the same in $W$ and
$W'$ except that $W$ traverses $2j$ more steps in color $a$, for some $j$.
Since $W$ and $W'$ have the same endpoints, this contradicts the lack of an
open parity walk under $f$.  Thus $W'$ does not exist.
\end{proof}

\section{Parity Edge-Coloring of Complete Graphs}

To prove that $\pecn(K_n)=2^{\CL{\lg n}}-1$ for all $n$, it suffices to prove
it when $n$ has the form $2^k+1$.  Below we prove it for $K_5$ and $K_9$ by
case analysis involving counting arguments.  Induced subgraphs of canonical
colorings provide the constructions; we only need the lower bounds.

\begin{proposition}\label{K5}
$\pecn(K_5)=7$.
\end{proposition}
\begin{proof}
Suppose that $K_5$ has a parity edge-coloring $f$ using at most six colors.
Each color class is a matching and hence has size at most 2.  Since $K_5$ has
10 edges, using at most six colors requires at least four color classes of
size 2.  Since any two colors used twice must not form a parity path of
length 4, each pair of colors used twice forms an alternating 4-cycle.
Hence the colors used twice are all restricted to the same four vertices.
However, there are only three disjoint matchings of size 2 in $K_4$.
Thus $f$ cannot exist.
\end{proof}

\begin{theorem}\label{K9}
$\pecn(K_9)=15$.
\end{theorem}
\begin{proof}
Let $f$ be a parity edge-coloring using at most 14 colors; we obtain a
contradiction.  Let $C_i$ be the set of edges in the $i$th color class,
and let $G_{i,j}$ be the spanning subgraph with edge set $C_i\cup C_j$.
By Lemma~\ref{connected}, a connected subgraph using any $k$ colors has at most
$2^k$ vertices.  Hence each $G_{i,j}$ has at least three components.
If $\C{C_i\cup C_j}\ge 7$, then $G_{i,j}$ has at most three components, since
the only non-tree components are 4-cycles, allowing the edges to be ordered
so that the first six edges reduce the number of components when added.

If each $G_{i,j}$ has at least four components, then
$\C{C_i\cup C_j}\le 6$.  If any class has size 4, then the others have size
at most two.  Since $K_9$ has 36 edges, and $4+2\cdot 13=30<36$, always
$\C{C_i}\le 3$.  Furthermore, since $7\cdot 3+7\cdot 2<36$, at least eight
classes have size 3; let $C_1$ be one of them.  If also $\C{C_j}=3$, then
$G_{1,j}$ has a 4-cycle, since otherwise six edges reduce $G_{1,j}$ to three
components.  The three edges of $C_1$ can form at most six 4-cycles with other
colors, but seven other classes have size 3.  The contradiction eliminates this
case.

Hence we may assume that $G_{1,2}$ has three components $A_1,A_2,A_3$ with
vertex sets $V_1,V_2,V_3$ and $\C{V_1}\le \C{V_2}\le \C{V_3}\le 4$.  Note that
$\C{V_2}\ge3$.  We show that for $i<j$, at least four colors join $V_i$ to
$V_j$.  If $\C{V_j}=4$, then the edges from $V_j$ to a vertex of $V_i$ have
distinct colors.  If $\C{V_j}<4$, then $\C{V_j}=3$ and $\C{V_i}\ge2$.  The
edges joining two vertices of $V_i$ to $V_j$ form $K_{2,3}$.  By
Corollary~\ref{K2n}, $\pecn(K_{2,3})=4$.

No color class $j$ outside $\{1,2\}$ connects one of $\{V_1,V_2,V_3\}$ to the
other two, since that would yield a connected 9-vertex graph in the three
colors $\{1,2,j\}$, contradicting Corollary~\ref{connected}.  With three
disjoint sets of four colors joining the pairs of components of $G_{1,2}$, we
now have 14 colors in $f$.  To avoid using another color, the remaining edges
joining vertices within components of $G_{1,2}$ must have colors used joining
those components.

Since $\C{V_2}\ge 3$, we may choose $u,v,w\in V_2$ with $uv\in C_1$ and
$vw\in C_2$ and $uw\in C_3$ and $e$ being an edge of $C_3$ that connects
distinct sets $V_i$ and $V_j$.
Suppose first that $e$ is incident to $V_2$.  If $\C{V_2}=4$, then $wx\in C_1$
or $ux\in C_2$, and appending $e$ to one end of $vu,uw,wx$ or $vw,wu,ux$ yields
a parity path.  If $\C{V_2}=3$, then $\C{V_1}\ge 2$, and the end of $e$ other
than $v$ is incident to an edge $e'$ in $C_1$ or $C_2$.  Now $e',e,vu,uw$ or
$e',e,vw,wu$ is a parity path.

Hence the endpoints of $e$ are in $V_1$ and $V_3$.  Let $z$ be the endpoint in
$V_3$.  If $\C{V_3}=4$, then each of the four colors joining $V_2$ to $V_3$
appears at each vertex of $V_2$.  Thus the color on $uz$ appears also on some
edge $wy$, and $e,zu,uw,wy$ is a parity path.

Hence $\C{V_1}=\C{V_2}=\C{V_3}=3$.  Since the nine edges joining $V_2$ and
$V_3$ use only four colors, some color is used on three of the edges.  Call it
$C_4$, with edges $uu',vv',ww'$ joining $V_2$ and $V_3$.  Avoiding a parity
path using $C_4$ with $C_1$ or $C_2$ forces $u'v'\in C_1$ and $v'w'\in C_2$.
If $z\in\{u',w'\}$, then $e,zu,uw,ww'$ or $e,zw,wu,uu'$ is a parity path.
Hence $z$ must be $v'$, and so $C_3$ appears only once on the copy of $K_{3,3}$
joining $V_2$ and $V_3$.

However, $K_{3,3}$ has no parity 4-edge-coloring in which one color is used
only once.  The other three colors would have multiplicities $3,3,2$.  Two
matchings of size 3 in $K_{3,3}$ form a 6-cycle, and a 2-colored 6-cycle
contains a parity path.
\end{proof}

It may be possible to generalize these arguments, but the case analysis
seems likely to grow.  Instead, we suggest another approach that could lead
to proving $\pecn(K_n)\ge 2^{\CL{\lg n}}-1$.

\begin{definition}\label{kset}
\rm
A {\it parity $r$-set edge-coloring} of a graph $G$ is an assignment of an
$r$-set of colors to each edge of $G$ so that every selection of a color
from the set on each edge yields a parity edge-coloring of $G$.
Let $\pecn_r(G)$ be the minimum size of the union of the color sets in
a parity $r$-set edge-coloring of $G$.
\end{definition}

Parity $r$-set edge-coloring is related to parity edge-coloring as $r$-set
coloring is to ordinary proper coloring.  An {\em $r$-set coloring} of a graph
assigns $r$-sets to the vertices so that the sets on adjacent vertices are
disjoint, with $\chi_r(G)$ being the least size of the union of the sets.  The
{\em $r$-set edge-chromatic number} $\chi_r'(G)$ is defined by
$\chi_r'(G)=\chi_r(L(G))$.  Thus $\pecn_r(G)\ge \chi_r'(G)$.

By using $r$ copies of an optimal parity edge-coloring with disjoint color
sets, it follows that $\pecn_r(G)\le r\pecn(G)$.  We have no examples yet
where equality does not hold.  Proving equality could help determine
$\pecn(K_n)$ by using the following result.

\begin{proposition}
If $K_n$ has an optimal parity edge-coloring in which some color class has
size $\FL{n/2}$, then $\pecn(K_n)\ge1+\pecn_2(K_{\CL{n/2}})$.
\end{proposition}
\begin{proof}
Let $f$ be an optimal parity edge-coloring in which $c$ is used on $\FL{n/2}$
edges.  Let $u_1v_1,\ldots,u_{\FL{n/2}}v_{\FL{n/2}}$ be the edges with color
$c$, and let $u_{\CL{n/2}}$ be the vertex missed by $c$ if $n$ is odd.
Contracting these edges yields $K_{\CL{n/2}}$, with $u_iv_i$ contracting to
$w_i$ for $i\le \FL{n/2}$, and $w_{\CL{n/2}}=u_{\CL{n/2}}$ when $n$ is odd.

Form a $2$-set edge-coloring $f'$ of $K_{\CL{n/2}}$ as follows.  For
$i<j\le \CL{n/2}$, let $f'(w_iw_j)=\{f(u_iu_j),f(v_iu_j)\}$.
Since $f'$ does not use $c$, to prove $\pecn_2(K_{\CL{n/2}})\le \pecn(K_n)-1$
it suffices to show that $f'$ is a parity $2$-set edge-coloring.

If $f'$ is not a parity $2$-set edge-coloring, then there is a parity
path $P'$ under some selection of edge colors from $f'$.  Form a path $P$ in
$K_n$ as follows.  When $P'$ follows the edge $w_iw_j$ with chosen color $a$,
$P$ moves along the edge $u_iv_i$ of color $c$ (if necessary) to reach an
endpoint in $\{u_i,v_i\}$ of an edge with color $a$ under $f$ whose other
endpoint is in $\{u_j,v_j\}$, and then it follows that edge.  This path has the
same usage as $P'$ for every color other than $c$.  Since $c$ misses only one
vertex of $K_n$, at least one end of $P'$ is a contracted vertex, and an edge
of color $c$ can be added or deleted at that end of $P$ to make the usage of
$c$ even if it had been odd.  If $P'$ is a $w_i,w_j$-path, then $P$ starts in
$\{u_i,v_i\}$ and ends in $\{u_j,v_j\}$ (one of the sets may degenerate to
$\{u_{\CL{n/2}}\}$).  Thus $P$ is a parity path under $f$, which is a
contradiction.
\end{proof}

If $n=2^k+1$, then $\CL{n/2}=2^{k-1}+1$.  If there is always an optimal
parity edge-coloring of $K_n$ with a near-perfect matching, then proving
$\pecn_2(K_n)=2\pecn(K_n)$ would inductively prove that
$\pecn(K_n)=2^{\CL{\lg n}}-1$.  Although we do not know whether
$\pecn_2(G)=2\pecn(G)$ in general, we provide support for the various
conjectures by proving this when $G$ is a path.

\begin{theorem}
$\pecn_r(P_n) = r\pecn(P_n)$.
\end{theorem}
\begin{proof}
We prove the stronger statement that for every parity $r$-set edge-coloring
$f$ of $P_n$, there is a set of $\pecn(P_n)$ edges whose color sets are
pairwise disjoint.

Let $\VEC e1{n-1}$ be the edge set of $P_n$ in order.  We say that a subset
$\{e_{i_1},\ldots,e_{i_q}\}$ of $E(P_n)$ with $\VECOP i1q<$ is {\em linked}
by $f$ if $f(e_{i_j})\cap f(e_{i_{j+1}})\ne\nul$ for $1\le j\le q-1$.

We claim first that if $E(P_n)$ decomposes into linked sets $\VEC S1t$ under
$f$, then setting $f'(e)=i$ when $e\in S_i$ yields a parity edge-coloring $f'$
of $P_n$ with $t$ colors.  If not, then there is a parity path $Q$ under $f'$.
Since $Q$ has an even number of edges with color $i$, we can pair successive
edges in the list of edges having color $i$ (first with second, third with
fourth, etc.).  Since $S_i$ is linked, we can pick a common color from the two
sets assigned to a pair.  Doing this for each pair and each color under $f'$
selects colors from the sets assigned to $Q$ under $f$ that form a parity path.
This contradicts the choice of $f$ as a parity $r$-set edge-coloring.  Thus
every partition of $E(P_n)$ into linked sets needs at least $\pecn(P_n)$ parts.

To obtain edges with disjoint color sets from such a partition, we first
construct a bipartite graph $H$ with partite sets $\VEC v1{n-1}$ and
$\VEC w1{n-1}$ by letting $v_iw_j$ be an edge if and only if $i < j$ and
$f(e_i) \cap f(e_j) \ne \nul$.  If $E(P_n)$ has a partition into $t$ linked
sets, then $H$ has a matching of size $n-1-t$, obtained by using the edge
$v_iw_j$ when $e_i$ and $e_j$ are successive elements in one part of the
partition.

The construction of a matching from a partition is reversible.  As edges are
added to the matching, starting from the empty matching and the partition into
singletons, the structural property is maintained that for the edges in a part,
only the first edge $e_j$ has $w_j$ unmatched, and only the last edge $e_i$ has
$v_i$ unmatched.  Hence when an edge $v_iw_j$ is added to the matching, it
links the end of one part to the beginning of another part, reduces the number
of parts, and maintains the structural property.

Thus $E(P_n)$ has a partition into $t$ linked sets under $f$ if and only if $H$
has a matching of size $n-1-t$.  When $t$ is minimized, the K\"onig--Egerv\'ary
Theorem yields a vertex cover of $H$ with size $n-t-1$.  Because the complement
of a vertex cover is an independent set, $H$ has an independent set $T$ of size
$n+t-1$.  Since $V(H)$ consists of $n-1$ pairs of the form $\{v_i,w_i\}$, at
least $t$ such pairs are contained in $T$.  If
$\{v_i,w_i\},\{v_j,w_j\}\esub T$, then $f(e_i)$ and $f(e_j)$ are disjoint.
Therefore there is a set of $t$ edges whose color sets are pairwise disjoint.

We conclude that $\pecn_r(P_n)\ge rt \ge r\pecn(P_n)$.
\end{proof}

\vspace{-1pc}

\section{Other Related Edge-Coloring Parameters}

In this section we describe other parameters defined by looser or more
restricted versions of parity edge-coloring, and we give examples to show
that $\pecn(G)$ is a different parameter.

A {\em nonrepetitive edge-coloring} is an edge-coloring in which no
pattern repeats immediately on a path.  That is, no path may have colors
$\VEC c1k,\VEC c1k$ in order for any $k$.  The notion was introduced for
graphs in \cite{AGHR}.  Every parity edge-coloring is nonrepetitive, and every
nonrepetitive edge-coloring is proper, so the minimum number of colors in a
non-repetitive edge-coloring of $G$ lies between $\pecn(G)$ and $\chi'(G)$.
The resulting parameter is called the {\em Thue chromatic number} in honor of
the famous theorem of Thue constructing non-repetitive sequences (generalized
to graphs in \cite{AGHR}).  The concept is surveyed in \cite{G1}.

More restricted versions of parity edge-colorings have also been studied.
A {\it conflict-free coloring} is an edge-coloring in which every path uses
some color exactly once.  An {\it edge-ranking} is an edge-coloring in which on
every path, the highest-indexed color appears exactly once.  Letting $c(G)$ and
$t(G)$ denote the minimum numbers of colors in a conflict-free coloring and an
edge-ranking, respectively, we have $t(G)\ge c(G)\ge \pecn(G)$.

Conflict-free coloring has been studied primarily in geometric settings; see
\cite{ELRS,HS,PT}.  Edge-rankings were introduced in \cite{IRV}.  It is known
that $t(K_n)\in\Omega(n^2)$ \cite{edgerank}; since $\pecn(K_n)\le 2n-3$, the
gap here can be large.  Equality can hold:
$t(P_n)=c(P_n)=\pecn(P_n)=\CL{\lg n}$.  Although computing $\pecn(G)$ or
$\specn(G)$ is NP-hard when $G$ is restricted to trees, there is a algorithm to
compute $t(G)$ that runs in linear time when $G$ is a tree \cite{LYtree} (at
least four slower polynomial-time algorithms were published earlier).
Computing $t(G)$ is NP-hard on general graphs \cite{LY}, as is finding a
spanning tree $T$ with minimal $t(T)$ \cite{MUI}.

In this string of inequalities, $c(G)$ and $\pecn(G)$ are neighboring
parameters.  In this section, we present examples to show that they may differ.
In fact, in all these examples $c(G)>\specn(G)$.

\begin{corollary}
$c(K_{2^k})>\specn(K_{2^k})=\pecn(K_{2^k})$ when $k\ge4$.
\end{corollary}

\vspace{-.5pc}

\begin{proof}
Let $f$ be an optimal spec of $K_{2^k}$.  By Theorem~\ref{canonical},
$f$ is canonical.  The spanning subgraph of $K_{2^k}$ formed by the color
classes whose names are vectors of weight 1 is isomorphic to the hypercube
$Q_k$, and the colors on it correspond to the coordinate directions.  If there
is a path in $Q_k$ that crosses each coordinate direction more than once, then
$f$ is not conflict-free and $c(K_{2^k})>\specn(K_{2^k})$.  In fact, it is
easy to find such paths when $k\ge4$.
\end{proof}

\begin{example}
\rm
As noted in Corollary~\ref{pathcyc}, $\specn(C_8)=3$.  Suppose that $C_8$
has a conflict-free 3-edge-coloring.  If a color is used only once, then the
other two colors alternate on paths of length 4 avoiding it, thus forming
parity paths of length 4.  Hence the sizes of the three color classes must be
$(4,2,2)$ or $(3,3,2)$.  Now deleting a edge from a largest color class yields
a spanning path on which no color appears only once.
\end{example}

By induction on the length, every path has an optimal parity edge-coloring that
is conflict-free (use a color only on a middle edge and apply the induction
hypothesis to each component obtained by deleting that edge).  This statement
does not hold for trees.

\begin{definition}\label{broom}
\rm
A {\em broom} is a tree formed by identifying an endpoint of a path with a
vertex of a star.  Let $T_k$ be the broom formed using $P_{2^k-2k+2}$ and a
leaf of a star with $k$ edges.  The {\em parity} of a vertex in $Q_k$ is the
parity of the weight of the $k$-tuple naming it.
\end{definition}

We prove that $T_k$ embeds in $Q_k$ but needs more than $k$ colors for a
conflict-free edge-coloring (for $k\ge4$).  W. Kinnersley (private
communication) showed this for $k=5$, and D. Cranston participated in early
discussions about the proofs.  We must first show that $T_k$ indeed
embeds in $Q_k$.  This follows from the embeddability (as spanning trees) of
equitable double-starlike graphs proved in \cite{KM}, since adding $k-2$ leaf
neighbors to the 2-valent neighbor of the $k$-valent vertex in $T_k$ yields
such a tree.  Their proof is lengthy; we give a short direct proof for this
special case.

\begin{lemma}\label{pathQ}
If $x$ and $y$ are distinct vertices of $Q_k$ having the same parity, then
there is a path of length $2^k-3$ in $Q_k$ that starts at $x$ and avoids $y$.
\end{lemma}
\begin{proof}
It is well known that $Q_k$ has a spanning cycle when $k\ge2$.  Since $Q_k$ is
edge-transitive, there is a spanning path from each vertex to any adjacent
vertex (for $k\ge 1$).

The desired path exists by inspection when $k=2$.  For larger $k$, we proceed
inductively.  Vertices $x$ and $y$ differ in an even number of bits; by
symmetry, we may assume that they differ in the first two bits.  Let $Q'$ and
$Q''$ be the $(k-1)$-dimensional subcubes induced by the vertices with first
bit 0 and first bit 1, respectively.  We may assume that $x\in V(Q')$.  There
is a spanning $x,u$-path $P'$ of $Q'$, where $u$ is the neighbor of $x$ obtained
by changing the third bit.  Note that $P'$ has length $2^{k-1}-1$.

Let $v$ be the neighbor of $u$ in $Q''$.  Since $v$ has the same parity as $y$,
and $v\ne y$, the induction hypothesis yields a path $P''$ of length
$2^{k-1}-3$ in $Q''$ that starts at $v$ and avoids $y$.  Together, $P'$, $uv$,
and $P''$ complete the desired path in $Q_k$.
\end{proof}

\begin{lemma}\label{broomQ}
For $k\ge2$, the broom $T_k$ embeds in $Q_k$, and hence
$\specn(T_k)=\pecn(T_k)=k$.
\end{lemma}
\begin{proof}
Note that $T_k=P_4\esub Q_k$ when $k=2$; we proceed inductively.
For $k>2$, the tree $T_k$ contains $T_{k-1}$, obtained by deleting one
leaf incident to the vertex $v$ of degree $k$ and $2^{k-1}-2$ vertices from the
other end.  With $Q'$ and $Q''$ defined in Lemma~\ref{pathQ}, by the induction
hypothesis $T_{k-1}$ embeds in $Q'$.  The distance in $T_{k-1}$ from $v$
to its leaf nonneighbor $u$ is $2^{k-1}-2(k-1)+2$.  This is even, so $u$ and
$v$ have the same parity.  Let $x$ and $y$ be the neighbors of $u$ and $v$
in $Q''$, respectively; also $x$ and $y$ have the same parity.  By
Lemma~\ref{pathQ}, $Q''$ contains a path $P$ of length $2^{k-1}-3$ starting
from $x$ and avoiding $y$.  Now adding $vy$, $ux$, and $P$ to the embedding
of $T_{k-1}$ yields the desired embedding of $T_k$ in $Q_k$.
\end{proof}

\begin{theorem}
If $k\ge4$, then $c(T_k)=k+1=\specn(T_k)+1$.
\end{theorem}
\begin{proof}
For $k=4$, a somewhat lengthy case analysis is needed to show $c(T_4)>4$;
we omit this.  Let $x$ be the vertex of degree $k$ in $T_k$.

For $k\ge5$, we decompose $T_k$ into several pieces.  At one end is a star $S$
with $k-1$ leaves and center $x$.  Let $P$ be the path of length $2^{k-2}$
beginning with $x$.  Let $R$ be the path of length $2^{k-1}$ beginning at the
other end of $P$.  Since $k\ge 5$, we have $2^{k-2}+2^{k-1}\le 2^k-2k+2$, so
$P$ and $R$ fit along the handle of the broom.  Ignore the rest of $T_k$ after
the end of $R$.

Consider a conflict-free $k$-edge-coloring of $S\cup P\cup R$.  Since $R$ has
$2^{k-1}+1$ vertices, at least $k$ colors appear on $E(R)$.  Since $P$ has
$2^{k-2}+1$ vertices, at least $k-1$ colors appear on $E(P)$.  Hence on
$P\cup R$ there are $k-1$ colors that appear at least twice, and only one color
$c$ appears exactly once.  Since $x$ has degree $k$, all $k$ colors appear
incident to $x$, including $c$.  Hence $c$ appears on some edge of $S$, and
adding this edge to $P\cup R$ yields a path on which every color appears at
least twice.

For $k\ge2$, we obtain a conflict-free $(k+1)$-edge-coloring using color $k+1$
only on the edge $e$ of $P$ at $x$.  Deleting $e$ leaves the star $S$ and a
path $P'$ with $2^k-2k+2$ vertices.  Since $S$ has $k-1$ edges and the length
of $P'$ is less than $2^k$, each has a conflict-free edge-coloring using colors
1 through $k$.  Paths from $V(S)$ to $V(P')$ use color $k+1$ exactly once.
\end{proof}

It remains unknown how large $c(G)$ can be when $\specn(G)=k$ or
$\pecn(G)=k$, either in general or when $G$ is restricted to be a tree.

\section{Open Problems}

Many interesting questions remain about parity edge-coloring and strong
parity edge-coloring.  We have already mentioned several and collect them
here with additional questions.

The first conjecture is known to be true for $n\le16$.  For the second,
we know the value of $\specn(K_{n,n})$ when $n\le8$.

\begin{conjecture}
$\pecn(K_n)=2^{\CL{\lg n}}-1$.
\end{conjecture}

\begin{conjecture}\label{Knn}
$\pecn(K_{n,n})=\specn(K_{n,n})=2^{\CL{\lg n}}$.
\end{conjecture}

We have found families of graphs $G$ such that $\specn(G)>\pecn(G)$
(see Example~\ref{phnep}), but they all contain odd cycles.

\begin{conjecture}\label{bipartite}
\rm
$\pecn(G)=\specn(G)$ for every bipartite graph $G$.
\end{conjecture}

\begin{question}
\rm
What is the maximum of $\specn(G)$ when $\pecn(G)=k$?
\end{question}

When $k=1$, we have $\specn(K_{k,n})=n$, which may be only a bit more than half
of $2^{\CL{\lg n}}$.  The growth in terms of $k$ is not yet known, although
Theorem~\ref{Kmn} sheds some light.

\begin{question}
\rm
Fix $n$.  For $k<n$, how do $\pecn(K_{k,n})$ and $\specn(K_{k,n})$ grow with
$k$?  In particular, if Conjectures~\ref{Knn} and \ref{bipartite} are true,
then what is $\min\{k\st \pecn(K_{k,n})=\specn(K_{k,n})=2^{\CL{\lg n}}\}$?
Does equality hold in Corollary~\ref{morebicliq}?
\end{question}

Using distinct colors yields $\pecn(T)\le n-1$ for every $n$-vertex tree, and
$K_{1,n-1}$ is the only $n$-vertex tree achieving equality.  More detailed
questions can be studied.  For example, a necessary condition for embedding a
tree $T$ in $Q_k$ is that $\Delta(T)\le k$, but this is not sufficient.

\begin{question}
\rm
What is the maximum of $\pecn(T)$ among $n$-vertex trees $T$ with $\Delta(T)=D$?
\end{question}

We observed from Theorem~\ref{tree} that testing $\pec(T)\le k$ is NP-hard.
This suggests complexity questions for more restricted problems.

\begin{question}
\rm
Do polynomial-time algorithms exist for computing $\pecn(T)$ on trees with
maximum degree $D$ or on trees with bounded diameter?
\end{question}

Theorem~\ref{checspec} shows that recognition of specs is in P.  However, we
do not know whether this holds for parity edge-coloring.  (It does hold for
edge-colorings of trees using the labeling procedure of Theorem~\ref{tree}.)

\begin{question}
\rm
What is the complexity of testing for a parity path in an edge-coloring?
\end{question}

Paths and complete graphs show that $\pecn(G)$ cannot be bounded by bounding
the maximum degree or the diameter.  However, bounding both parameters
limits the number of vertices.  Hence the next question makes sense.

\begin{question}
\rm
What is the maximum of $\pecn(G)$ among graphs (or trees) with
$\Delta(G)\le k$ and $\diam(G)\le d$?
\end{question}

It is a classical question to determine the maximum number of edges in
an $n$-vertex subgraph of $Q_k$, where $n\le 2^k$.  Does the resulting graph
have the maximum number of edges in an $n$-vertex graph with parity
edge-chromatic number $k$?  More generally,

\begin{question}
\rm
What is the minimum of $\pecn(G)$ among all $n$-vertex graphs $G$
having minimum degree $t$?  Among those having $m$ edges?
\end{question}

The lower bound in Corollary~\ref{connected} naturally leads us to ask which
graphs achieve equality.  Every spanning subgraph of a hypercube satisfies
$\pecn(G)=\lg n(G)$; is the converse true?

\begin{question}
\rm
Which connected graphs $G$ satisfy $\pecn(G) = \CL{\lg n(G)}$?  Which satisfy
$\specn(G) = \CL{\lg n(G)}$?
\end{question}

Motivated by the uniqueness of the optimal spec of $K_{2^k}$, Dhruv Mubayi
suggested studying the ``stability'' of the result.

\begin{question}
\rm
Does there exist an parity edge-coloring of $K_{2^k}$ with $(1+o(1))2^k$
colors that is ``far'' from the canonical coloring?
\end{question}

In Section 5, we showed that paths satisfy all three properties below.
Are there other such graphs?

\begin{question}
\rm
For which graphs $G$ do the following (successively stronger) properties hold?
(a) $\pecn_2(G)=2\pecn(G)$?\\
(b) $\pecn_r(G)=r\pecn(G)$ for all $r$?\\
(c) every parity $r$-set edge-coloring of $G$ contains a set of $\pecn(G)$
edges whose color sets are pairwise disjoint?
\end{question}

Lemma~\ref{addedge}(a) states that deleting an edge reduces the parity
edge-chromatic number by at most 1.  Ordinary coloring has the same property.
Thus we are motivated to call a graph $G$ {\em critical} if
$\pecn(G-e)<\pecn(G)$ for all $e\in E(G)$.  We say that $G$ is {\em
doubly-critical} if $\pecn(G-e-e')=\pecn(G)-2$ for all $e,e'\in E(G)$.  Our
results on paths and cycles imply that for all $n \geq 1$, $P_{2^n+1}$ is
critical and $C_{2^n+1}$ is doubly-critical.  Naturally, any star is
doubly-critical.

\begin{question}
\rm
Which graphs are critical?  Which graphs are doubly-critical?
\end{question}

Since the factors can be treated independently in constructing a spec,
$\specn$ is subadditive under Cartesian product.
Note that $\specn(P_2\Box P_2)=2=\specn(P_2)+\specn(P_2)$.

\begin{question}
\rm
For what graphs $G$ and $H$ does equality hold in
$\specn(G\Box H)\le \specn(G)+\specn(H)$?  What can be said about
$\pecn(G\Box H)$ in terms of $\pecn(G)$ and $\pecn(H)$?
\end{question}

It may be interesting to compare $\pec(G)$ with related parameters
such as conflict-free edge-chromatic number on special classes of graphs.
We suggest two specific questions.

\begin{question}
\rm
What is the maximum of $c(T)$ such that $T$ is a tree with $\pecn(T)=k$?
What is the maximum among all graphs with parity edge-chromatic number $k$?
\end{question}

Finally, the definitions of parity edge-coloring and spec extend naturally to
directed graphs: the parity condition is the same but is required only for
directed paths or walks.  Hence $\pecn(D)\le \pecn(G)$ and
$\specn(D)\le \specn(G)$ when $D$ is an orientation of $G$.

For a directed path $\vec P_m$, the constraints are the same as for an
undirected path.  More generally, if $D$ is an acyclic digraph, and $m$ is the
maximum number of vertices in a path in $D$, then
$\pecn(D)=\specn(D)=\CL{\lg m}$.  The lower bound is from any longest path.

For the upper bound, give each vertex $x$ a label $l(x)$ that is the maximum
number of vertices in a path ending at $x$ (sources have label $0$).
Write each label as a binary $\CL{\lg m}$-tuple.  By construction,
$l(v)>l(u)$ whenever $uv$ is an edge.  To form a spec of $D$, use a color $c_i$
on edge $uv$ if the $i$th bit is the first bit where $l(u)$ and $l(v)$ differ.
All walks are paths.  Any $x,y$-path has odd usage of $c_i$, where the $i$th is
the first bit where $l(x)$ and $l(y)$ differ, since no edge along the path
can change an earlier bit.

Thus the parameters equal $\CL{\lg n}$ for the $n$-vertex transitive
tournament, which contains $\vec P_n$.  This suggests our final question.

\begin{question}
\rm
What is the maximum of $\pecn(T)$ or $\specn(T)$ when $T$ is an
$n$-vertex tournament?
\end{question}

\bigskip
\centerline{\Large Acknowledgement}
\medskip
We thank Dan Cranston, Will Kinnersley, Brighten Godfrey, Michael Barrus, and
Mohit Kumbhat for helpful discussions.

\end{document}